\documentclass{article}

\usepackage{amsfonts}
\usepackage[T1]{fontenc}
\usepackage{array}
\usepackage{tabularx}
\usepackage{amssymb}
\usepackage{tikz-cd}

\usepackage{graphicx}
\def\vnabla{^{\vup}\!\!\nabla}
\def\vup{^{\vee}\!\!\!}

\def\C{\mathbb{C}} 
\def\G{\mathbb{G}} 

\def\R{\mathbb{R}}  

\def\no{\noindent}

\def\beq{\begin{equation}}
\def\eeq{\end{equation}}

\def\w{\wedge}

\begin{document}

\title{Geometric Algebras of Light Cone Projective Graph Geometries}

\author{Garret Sobczyk\\ Departamento de Actuaría F\'isica y Matem\'aticas,\\Universidad de las Am\'ericas-Puebla,\\72820 Puebla, Pue., M\'exico\\
	garretudla@gmail.com\\
}
\maketitle              
 
\begin{abstract}

A null vector $v$ is an algebraic quantity with the property that $v^2=0$. I denote the universal algebra generated by taking all sums and products of null vectors over the real or complex numbers  by ${\cal N}$, \cite{GS2022}. The rules of addition and multiplication in ${\cal N}$ are taken to be the familiar rules of addition and multiplication of real or complex square matrices.
 A pair of null vectors is {\it positively} or {\it negatively} correlated if their inner product is {\it positive} or {\it negative}, respectively.
A large class of geometric algebras are isomorphic to real or complex $2^n \times 2^n$ matrix algebras, or a pair of such algebra. 
I begin the study of the eigenvector-eigenvalue problem of linear operators in the geometric algebra $\G_{1,n}$ of $\R^{n+1}$, and by restricting to barycentric coordinates, $n$-{\it simplices} whose $n+1$ vertices are non-zero null vectors.
  These ideas provide a foundation for a new {\it Cayley-Grassmann Theory of Linear Algebra}, with many possible applications in pure-applied areas of science and engineering. 
 
 	\smallskip
 \no {\em AMS Subject Classification:} 03B30,05C20,15A66,15A75
 \smallskip
 
 \no {\em Keywords:} Clifford algebra, complete graphs, Grassmann algebra, Lorentzian spacetime.

\end{abstract}

\section*{0. Introduction}

The origin of the ideas in this paper date back most directly to mathematics that was set down in the nineteenth century by H. Grassmann \cite{Dieu}, A. Cayley {\it (Memoir on the
Theory of Matrices 1858)}, and W. Clifford \cite{WKing}. It is regrettable today, after more than $150$ years, that Clifford's {\it geometric algebra} has not found it proper place in the {\it Halls of Mathematics and Science} \cite{ICAtalk}. My journey in this saga began in 1965, when I starting working in geometric algebra as a graduate Ph.D. student of Professor David Hestenes at Arizona State University \cite{HS1992}, and continued with years spent with gracious colleagues in Poland and Mexico. It is my belief that this paper will bring us closer to the day when geometric algebra has finally found its proper place in the Millennial Human Quest for the development of the {\it geometric concept of number} \cite{SNF}.

In Section 1, it is shown that the geometric algebras
$\G_{1,n}$ and $\G_{n,1}$ have special bases of  all positively, or all negatively correlated null vectors, respectively. In the case of $\G_{1,n}$, the inner products can all be chosen to be $+\frac{1}{2}$, and in the case of $\G_{n,1}$, $-\frac{1}{2}$. For simplicity, the classification of endomorphisms on $\R^{n+1}$ is considered only in the case of a $(+\frac{1}{2})$-positively correlated basis of a geometric algebra $\G_{1,n}$, but the same analysis is valid for studying endomorphisms on $\R^{n+1}$ of a $(-\frac{1}{2})$-negatively correlated basis of a geometric algebra $\G_{n,1}$.  

   In Section 2, basic ideas of linear algebra in $\R^{n+1}$ are developed in the symmetric algebra ${\cal A}_{n+1}^+$ of a correlated basis of null vectors in $\G_{1,n}^1$. The concept of a LPGG {\it star projection} of a geometric number is defined and studied. The {\it vector derivative} $\nabla$ is defined, paying particular attention to its important properties.   
   
   In Section 3, basic properties of lower dimensional geometric algebras are explored in the correlated basis algebra ${\cal A}_3^+$ of $\G_{1,3}$. The concept of the LPGG star projection suggests that a new classification all geometric algebras is possible in the correlated null vector algebra
   ${\cal A}_{n+1}^+$ of $\G_{1,n}$, \cite{GS3}. 

In Section 4, by introducing barycentric coordinates, complete graphs are studied in which every pair of vertices is connected by an edge. {\it Light Cone Projective Geometry} (LPGG) 
 is built upon the property that for any dimension $n\ge 1$, there exits {\it positively}, or {\it negatively} {\it correlated light cones}, defined by sets of $(n+1)$ null {\it basis vectors} $\{a_1, \ldots, a_{n+1}\}$ of $\G_{1,n}$, or $\G_{n,1}$, such that $a_1\w \cdots \w a_{n+1} \ne 0$ and $a_i\cdot a_j =\pm\frac{(1-\delta_{ij})}{2}$, respectively.

\section{The geometric algebras $\G_{1,n}$ and $\G_{n,1}$ of $\R^{n+1}$}
The geometric algebras
$\G_{1,n}$ and $\G_{n,1}$ arise from null vector bases of $ \R^{n+1}$ by constructing positively, or negatively correlated, null vectors in terms of the {\it standard bases} $\{e_1,f_1, \cdots, f_n\}$ of $\G_{1,n}$, or $\{f_1,e_1, \cdots, e_n\}$ of $\G_{1,n}$, respectively. Renewed interest in these Clifford algebras is due in part to the pivotal Lecture Notes published by Marcel Riesz in 1958, \cite{mr1958}.   
The geometric algebras $\G_{1,n}$ and $\G_{n,1}$ make up the two {\it fundamental sequences} 
of successively larger algebras, 
\beq \R \subset \G_{1,1}  \subset \G_{1,2} \subset \G_{1,3} \subset \cdots \subset \G_{1,n} \subset \cdots \subset {\cal N},  \label{seq-GA1} \eeq
and 
\beq \R \subset \G_{1,1}  \subset \G_{2,1} \subset \G_{3,1} \subset \cdots \subset \G_{n,1} \subset \cdots \subset {\cal N},  \label{seq-GA2} \eeq 
where $\cal N$ is the universal algebra generated by taking sums and products of null vectors. 
  See \cite{GS2022,SNF,GS2}, and other references.

Let $\{a_1, \ldots a_{n+1}\}\subset \R^{n+1} $ be a set of positively, or negatively, correlated null vectors satisfying $a_1\w \cdots \w a_{n+1}\ne 0$, and the $(n+1)^2$ properties,
\beq a_i \cdot a_j\equiv \frac{1}{2}(a_i a_j + a_j a_i) :=\pm \frac{1-\delta_{ij}}{2} \quad {\rm for} \quad 1\le i,j \le n+1, \label{defprop} \eeq 
respectively, where $\delta_{ij}$ is the usual delta function. 
In terms of these basis null vectors, 
\beq  \R^{n+1} := \{x| \ x=x_1 a_1 + \cdots + x_{n+1} a_{n+1}, \ x_i\in \R  \}.\label{defRn} \eeq
 The multiplication tables for sets of positively (PC), or negatively (NC), correlated null vectors $a_i , a_j$, for $1 \le i<j \le n+1$, follow directly from the properties (\ref{defprop}), and generate the positively, and negatively correlated null vector algebras ${\cal A}_{1,n}^+=\G_{1,n} $, and ${\cal A}_{n,1}^-=\G_{n,1}$, respectively.
\begin{table}[h!]
	\begin{center}
		\caption{Multiplication Table.}
		\label{table3.1}
		\begin{tabular}{c|c|c|c|c} 
			& $a_i $  & $ a_j $  & $a_i a_j $  & $ a_j a_i $ \cr 
			\hline
			$ a_i$  &   0   & $ a_i a_j $     & $  0 $         & $ a_i $ \cr 
			\hline         
			$ a_j $  &  $ a_j a_i $  & $ 0 $     & $  a_j $         & $ 0 $ \cr 
			\hline     
			$ a_i a_j $    &  $a_i $    &  0  &  $a_i a_j $    & 0  \cr
			\hline
			$ a_j a_i $  &  0   &  $a_j $   &  0     & $a_j a_i$   \cr
		\end{tabular}
	\end{center}
\end{table}

For a set of positively or negatively correlated null vectors $\{a_1, \ldots, a_n\}$, define
\beq A_k := \sum_{i=1}^k a_i.  \label{sumak} \eeq
The geometric algebra 
\[ \G_{1,n}:=\R(e_1,f_1,\ldots, f_n) ,\] 
where $\{ e_1, f_1, \ldots, f_n \}$ is the {\it standard basis} of anticommuting orthonormal vectors, with $e_1^2=1$ and $f_1^2 = \cdots = f_n^2 =-1$.
The $2^{n+1}$-{\it canonical forms} of the standard multivector basis elements are
\beq \Big\{1; e_1,f_1,\cdots, f_n;e_1f_1,\cdots, e_1f_{n}, \big[_{1\le i<k\le n}f_if_k,\big]; \ldots;e_1f_1\cdots f_n\Big\}.   \label{Fcanform}\eeq

Alternatively, the geometric algebra $\G_{1,n}$ can be defined by
\[ \G_{1,n}:=\R(a_1,\ldots, a_{n+1})=:{\cal A}_{1,n}^+,   \]
where   
$\{a_1, \ldots, a_{n+1} \}$ is a set of positively correlated null vectors satisfying the Multiplication Table 1. In this case, the standard basis vectors of $\G_{1,n}$ can be defined by $e_1=a_1+a_2=A_2$, $f_1=a_1-a_2=A_1- a_2 $, and for $2\le k\le n$
\beq f_k=\alpha_k \Big(A_k-(k-1)a_{k+1}\Big),   \label{formulafk}  \eeq  
where $\alpha_k := \frac{-\sqrt 2}{\sqrt{k(k-1)}}$.
 The $2^{n+1}$-{\it canonical forms} of the standard multivector basis elements of ${\cal A}^+_{n+1}$ are
 \beq \Big\{1;a_1 , \ldots, a_{n+1}; \big[_{1\le i < j \le n+1}a_i a_j,\big];\ldots;a_1\cdots a_{n+1}\Big\}.   \label{Acanform}\eeq

The geometric algebra 
\[ \G_{n,1}:=\R(f_1,e_1,\ldots, e_n) ,\] 
where $\{ f_1, e_1, \ldots e_n \}$ is the {\it standard basis} of anticommuting orthonormal vectors, with $f_1^2=-1$ and $e_1^2 = \cdots = e_n^2 =1$. 
Alternatively, the geometric algebra $\G_{n,1}$ can be defined by
\[  \G_{n,1}:=\R(a_1,\ldots, a_{n+1})=:{\cal A}^-_{1,n},   \]
where   
$\{a_1, \ldots, a_{n+1} \}$ is a set of negatively correlated null vectors satisfying the Multiplication Table 2. In this case, the standard basis vectors of $\G_{n,1}$ can be defined by $f_1=a_1+a_2=A_2$, $e_1=a_1-a_2=A_1- a_2 $, and for $2\le k\le n$
\beq e_k=\alpha_k \Big(A_k-(k-1)a_{k+1}\Big), \label{formulaek} \eeq   
where $\alpha_k := \frac{-\sqrt 2}{\sqrt{k(k-1)}}$.
The $2^{n+1}$-{\it canonical forms} of the standard multivector basis elements of ${\cal A}^-_{n+1}$ is the same as (\ref{Acanform}).

\begin{table}[h!]
	\begin{center}
		\caption{Multiplication Table.}
		\label{table3.1}
		\begin{tabular}{c|c|c|c|c} 
			& $a_i $  & $ a_j $  & $a_i a_j $  & $ a_j a_i $ \cr 
			\hline
			$ a_i$  &   0   & $ a_i a_j $     & $  0 $         & $ -a_i $ \cr 
			\hline         
			$ a_j $  &  $ a_j a_i $  & $ 0 $     & $  -a_j $         & $ 0 $ \cr 
			\hline     
			$ a_i a_j $    &  $-a_i $    &  0  &  $-a_i a_j $    & 0  \cr
			\hline
			$ a_j a_i $  &  0   &  $-a_j $   &  0     & $-a_j a_i$   \cr
		\end{tabular}
	\end{center}
\end{table} 

For the remainder of this paper, only properties of the positively correlated null vector algebras $ {\cal A}^+_{n+1}:={\cal A}_{1,n}$ of the geometric algebras  of $\G_{1,n}$ are considered. It should be recognized, however, that any of these properties can be easily translated to the corresponding properties of the negatively correlated null vector basis algebras ${\cal A}^-_{n+1}:={\cal A}_{n,1}$ of $\G_{n,1}$. Indeed, much more general algebras of correlated null vectors in ${\cal N}$ can be defined and studied, but with correspondingly more complicated rules of multiplication. In addition to providing a new framework for the study of Linear Algebra on $\R^{n+1}$, the last section of the paper shows how the ideas can be applied to graph theory. 

\section{Linear algebra of $\R^{n+1}$ in ${\cal A}_{n+1}^+$}

The {\it position vector} $x \in \R^{n+1}$ in the standard basis of $\G_{1,n}$ is
\beq  x := s_1 e_1 + \sum_{i=1}^n s_{i+1} f_i \in \G_{1,n}^1 .\label{xinG1n} \eeq 
 Alternatively, in the correlated null vector basis algebra ${\cal A}^+_{n+1}=\G_{1,n}$,
 \beq x = \sum_{i=1}^{n+1} x_i a_i \in {\cal A}^+_{n+1}. \label{xinA1n} \eeq
 
Since geometric algebras are fully compatible with matrix algebras, matrix algebras over geometric algebras are well defined \cite{SNF}.
To relate the bases (\ref{xinG1n}) and (\ref{xinA1n}), in matrix notation
 \beq x = \pmatrix{s_1 & \ldots & s_{n+1}}
 \pmatrix{e_1 \cr f_1 \cr \cdot \cr \cdot \cr \cdot \cr f_n} = \pmatrix{x_1 & \ldots & x_{n+1}}
 \pmatrix{a_1 \cr a_2 \cr \cdot \cr \cdot \cr \cdot \cr a_{n+1}}, \label{positionx} \eeq
 or in {\it abbreviated form}, $x=s_{(n+1)}F_{(n+1)}=x_{(n+1)}A_{(n+1)}$.
 The {\it quadratic form} of $\G_{1,n}$ is specified by $F^{-1}_{(n+1)}:=F_{(n+1)}^tB$, where
 \beq B := F_{(n+1)}\cdot F_{(n+1)}^t =\pmatrix{e_1 \cr f_1 \cr \cdot \cr \cdot \cr \cdot \cr f_n} \cdot \pmatrix{e_1 & f_1 & \ldots & f_n}  \label{QformB} \eeq
 \[ = \pmatrix{1 & 0 & 0 & \cdots & 0 \cr 
               0 & -1 & 0 & \cdots & 0 \cr
          \cdot  &  \cdot & \cdot & \cdot & \cdot  \cr
           \cdot   &  \cdot & \cdot & \cdot & \cdot \cr
              0 & 0 & \cdots & 0 & -1} ,  \]
  where  $F_{(n+1)}^t$ denotes the {\it row transpose} of the column $F_{(n+1)}$.          
          
    Let $v,w\in \G_{1,n}^1$ be vectors. Expressed in the {\it standard basis} of $\G_{1,n}$, the geometric product
    \[  vw=v_{(n+1)}F_{(n+1)}F_{(n+1)}^tw_{(n+1)}^t
     \]
  \[  = v_{(n+1)}F_{(n+1)}\cdot F_{(n+1)}^tw_{(n+1)}^t+ v_{(n+1)}F_{(n+1)}\w F_{(n+1)}^tw_{(n+1)}^t \] 
  \beq =v_{(n+1)}B w_{(n+1)}^t+ v_{(n+1)} \pmatrix{0 & e_1f_1 & e_1f_2 & \cdots & e_1f_n \cr 
  	f_1e_1 & 0 &	f_1f_2 & \cdots & 	f_1f_n \cr
  	\cdot  &  \cdot & \cdot & \cdot & \cdot  \cr
  	\cdot   &  \cdot & \cdot & \cdot & \cdot \cr
  	f_ne_1 & f_n f_1 &	f_nf_2 & \cdots & 0 }w_{(n+1)}^t .     \label{vwgaprod} \eeq      
                   
Dotting each side of the equation (\ref{positionx}) on the right by the row matrix 
\[ F^{-1}_{(n+1)}:=\pmatrix{e_1 & -f_1 & \ldots & -f_n},\]
 and noting that 
 \[ \pmatrix{e_1 \cr f_1 \cr \cdot \cr \cdot \cr \cdot \cr f_n} \cdot \pmatrix{e_1 & -f_1 & \ldots & -f_n}   \]
is an expression for the $(n+1)\times (n+1)$ identity matrix, immediately gives $s_{(n+1)} = x_{(n+1)}T$, where the matrix of transition $T$ is defined by the {\it Gramian matrix} 
\beq T:=A_{(n+1)}\cdot F^{-1}_{(n+1)}=\pmatrix{a_1 \cr a_2 \cr \cdot \cr \cdot \cr \cdot \cr a_{n+1}}\cdot \pmatrix{e_1 & -f_1 & \ldots & -f_n}  \label{transmatrix} \eeq
 in terms of the inner products $a_i \cdot f_j$.  
 These inner products are directly calculated using (\ref{formulafk}). The transition matrix $T_8$, and its inverse $T_8^{-1}$, for the geometric algebra $\G_{1,7}$ is given in Appendix A.
 
Note, that whereas $T$ is the transition matrix
\beq TF_{(n+1)}=T \pmatrix{e_1 \cr f_1 \cr \cdot \cr \cdot \cr f_n }=\pmatrix{a_1 \cr \cdot \cr \cdot \cr a_{n+1}}=A_{(n+1)} , \label{etoa} \eeq
of the column basis vectors $F_{(n+1)}$ to the column basis vectors $A_{(n+1)}$,
$T^{-1}$ is the coordinate transition matrix
\beq s_{(n+1)}T^{-1}=\pmatrix{s_1 & \cdots & s_{n+1}}T^{-1}=\pmatrix{x_1 & \cdots & x_{n+1}} = x_{(n+1)}, \label{coorstox} \eeq 
from the {\it row vector} coordinates of $x$ to the
row vector coordinates of $x$ in the basis $ A_{(n+1)}$. Great care must be taken to avoid confusion.

Converting the calculation in (\ref{vwgaprod}) to a calculation for 
$v,w \in {\cal A}_{n+1}^+$,
 \[  vw=v_{(n+1)}F_{(n+1)}F_{(n+1)}^tw_{(n+1)}^t
\]
 \[ =v_{(n+1)}T^{-1}TF_{(n+1)}F_{(n+1)}^tT^t(T^t)^{-1}
  w_{(n+1)}^t
\]
\[ =v_{(n+1)}^aA_{(n+1)}A_{(n+1)}^t
(w_{(n+1)}^a)^t,
\]
giving
\[  vw = v_{(n+1)}^aA_{(n+1)}\cdot A_{(n+1)}^t(w_{(n+1)}^a)^t+ v_{(n+1)}^aA_{(n+1)}\w A_{(n+1)}^t(w_{(n+1)}^a)^t \]
 \beq = v_{(n+1)}^a \pmatrix{0 & a_1 a_2 & a_1 a_3 & \cdots & a_1 a_{n+1} \cr 
	a_2 a_1 & 0 &	a_2 a_3 & \cdots & 	a_2 a_{n+1} \cr
	\cdot  &  \cdot & \cdot & \cdot & \cdot  \cr
	\cdot   &  \cdot & \cdot & \cdot & \cdot \cr
	a_{n+1}a_1 & a_{n+1} a_2 &	a_{n+1}a_3 & \cdots & 0 }(w_{(n+1)}^a)^t .     \label{vwaaprod} \eeq      
 
 \subsection{Bivector endomorphisms in ${\cal  A}_{n+1}^+$}
 The standard treatment of the relationship between Clifford's geometric algebras and Cayley's matrix algebras is well-known, \cite[p.74]{GS2}, \cite[p.217]{PL2001}. Something that has always disturbed me is that this relationship is an isomorphism only for square matrix algebras of order $2^n\times 2^n$. This sorry state of affairs is at least partially rectified in the algebras ${\cal A}_{n+1}^+ = \G_{1,n}$.
 
 Recalling the definition of (\ref{sumak}), it is not difficult to show that for $k\ge 2$,
 \beq (A_k)^2=A_{k-1}^2+2\Big(\sum_{j=2}^{k}A_{j-1}^2\cdot a_j\Big)=(k-1)+\cdots +1 =\pmatrix{k \cr 2}.  \label{sqrAk}\eeq 
 Defining $\hat A_k :={\frac{\sqrt 2}{\sqrt{k(k-1)}}}A_k$, it follows that  $\hat A_k^2=1$. For $g\in {\cal A}_{n+1}^+$, and $k>1$, I now define the
 {\it LPGG k-projection} of $g$,
 \beq g^\star := \hat A_{k}g\hat A_{k}. \label{gstar} \eeq
  When $k=n+1$, the $\star$-projection becomes a {\it conjugation} on $\G_{1,n}$. 
   
 Clearly, for all $g\in {\cal A}_{n+1}^+$, $(g^\star)^\star = g$, and given a second $h\in  {\cal A}_{n+1}^+$, 
 \beq (gh)^\star = \hat A_{n+1}g\hat A_{n+1} \hat A_{n+1}h \hat  A_{n+1} =g^\star h^\star .               \label{ghaaprod} \eeq 
Defining the {\it A-matrix} of $g\in {\cal A}_{n+1}^+$,
\beq  [g]_a := A_{(n+1)}gA_{(n+1)}^t=[a_i g a_j]_a,  \label{Amatrixg} \eeq
it follows that
\[ g^\star = \hat A_{n+1} g {\hat A}_{n+1} = 
  {\cal I}_{(n+1)}^t [g]_a {\cal I}_{(n+1)},  \]
  where ${\cal I}_{(n+1)}$ and ${\cal I}_{(n+1)}^t$ are the $n+1$ column and row matrices 
  \[ {\cal I}_{(n+1)}:= \pmatrix{1 \cr \cdot \cr \cdot \cr \cdot \cr 1}, \ \ {\cal I}_{(n+1)}^t:= \pmatrix{1 & 1 &  \cdots & 1},  \] 
respectively.

The GA product of $g,h \in {\cal A}_{n+1}^+$, in terms of their $(n+1)\times (n+1)$ $A$-matrices $[g]_a,[h]_a$, then takes the unusual form
  \beq gh = \hat A_{n+1} {\cal I}_{(n+1)}^t \Big([g]_a {\cal I}_{(n+1)} {\cal I}_{(n+1)}^t [h]_a\Big) {\cal I}_{(n+1)} \hat A_{n+1},  \label{Aghproduct} \eeq
  mediated by the square singular $(n+1)$-matrix
  ${\cal I}_{(n+1)} {\cal I}_{(n+1)}^t $. Equation (\ref{Aghproduct}) is a generalization of the closely related formula (\ref{vwaaprod}) for the multiplication of the vectors $v,w \in {\cal A}_{n+1}^+$.
 
 It follows from (\ref{gstar}) and (\ref{Aghproduct}) that a real or complex $(n+1)\times (n+1)$-matrix $[g_{ij}]$ is the matrix of a scalar plus a bivector $g\in {\cal A}^+_{n+1}$, that is
 \[ [g]_a=\hat A_{n+1}[g^\star ]_a\hat A_{n+1} \]
  \beq =\pmatrix{0 & g_{12}a_1 a_2 &g_{13} a_1 a_3 & \cdots &g_{1,n+1} a_1 a_{n+1} \cr 
  	g_{21}a_2 a_1 & 0 &g_{23}	a_2 a_3 & \cdots & 	g_{2,n+1}a_2 a_{n+1} \cr
  	\cdot  &  \cdot & \cdot & \cdot & \cdot  \cr
  	\cdot   &  \cdot & \cdot & \cdot & \cdot \cr
  	g_{1,n+1}a_{n+1}a_1 & g_{2,n+1}a_{n+1} a_2 &	g_{n,n+1}a_{n+1}a_3 & \cdots & 0 } .  \label{rmatrixofg} \eeq
  Comparing the matrix $[g]_a$ in (\ref{rmatrixofg}) to the matrix in (\ref{vwgaprod}), seems to contradict  that the trace of a matrix is invariant under a change of basis. However this is not the case since the terms $g_{ij}a_ia_j$ of (\ref{rmatrixofg}) consists of scalars and bivectors. Pl\"ucker relations are important in understanding the structure of bivectors \cite{GS1}, particularly bivectors in $\G_{1,n}$, and in study of conformal mappings \cite{notices}.

  \subsection{The gradient $\nabla$ } 
   
  A crucial tool for carrying out calculations in the geometric algebra $\G_{1,n}$ is the {\it gradient} $\nabla$. In the references \cite{HS1992,SNF,HDVA1993}, the gradient $\nabla$, alongside the geometric algebra $\G_n$, has been developed as a basic tool for formulating and proving basic theorems of linear algebra in $\R^{n}$. Since the properties of the gradient are independent of the quadratic form of the geometric algebra used, instead of using the Euclidean geometric algebra $\G_{n+1}$ of $\R^{n+1}$, we can equally well define it in terms of the geometric algebra $\G_{1,n}$. It follows that all theorems of linear algebra developed in  \cite{HS1992,SNF,HDVA1993} are equally valid in $\G_{1,n}$ without modification. In the standard basis of $\G_{1,n}$,
  \beq \nabla := e_1 \frac{\partial}{\partial s_1} - f_1 \frac{\partial}{\partial s_2} - \cdots - f_n \frac{\partial}{\partial s_{n+1}}. \label{defgradient} \eeq 
  
  With the transition matrix (\ref{transmatrix}) in hand, the expression for the gradient in the null vector basis of ${\cal A}_{1,n}$,
  \beq \nabla = \sum_{i=1}^{n+1} (\nabla x_i)\frac{\partial}{\partial x_i}, \label{gradientx} \eeq 
is nothing more than a simple expression of the {\it chain rule} in calculus. In terms of the abbreviated notation for (\ref{positionx}), it is not difficult to derive the transformation rules relating the {\it bases columns}  $A_{(n+1)}$ and $F_{(n+1)}$. 

Using (\ref{transmatrix}), and solving
\beq x=s_{(n+1)}F_{(n+1)} = x_{(n+1)}A_{(n+1)},	\label{basesrel} \eeq 
gives the important relations
\begin{itemize}
   \item $s_{(n+1)}= x_{(n+1)}T \ \ \iff \ \   x_{(n+1)}=s_{(n+1)}T^{-1}$
	\item $x\cdot F_{(n+1)}^{-1}=s_{(n+1)}= x_{(n+1)}T \ \ \iff \ \ F_{(n+1)}^{-1} =\nabla s_{(n+1)}= \nabla x_{(n+1)}T$	
	\item $ x_{(n+1)}=s_{(n+1)}T^{-1} \ \ \iff \ \  A_{(n+1)}^{-1}:=\nabla x_{(n+1)}=F_{(n+1)}^{-1}T^{-1}$.
	\item $A_{(n+1)}^{-1}=F_{(n+1)}^{-1}T^{-1} \ \ \iff \ \ A_{(n+1)}^{-1}T=F_{(n+1)}^{-1}$.
	\item $\nabla x = n+1 = F_{(n+1)}^{-1} F_{(n+1)} = A_{(n+1)}^{-1}A_{(n+1)}$. 
    \item $F_{(n+1)}=F_{(n+1)} \cdot \nabla x =F_{(n+1)}\cdot \nabla x_{(n+1)}A_{(n+1)}$
    
    \quad \quad\quad $=F_{(n+1)}\cdot \nabla s_{(n+1)}T^{-1}A_{(n+1)} = T^{-1} A_{(n+1)}$.
\end{itemize}

Note, whereas $F_{(n+1)}$ and $A_{(n+1)}$ have been defined as {\it column} matrices of vectors, $F_{(n+1)}^{-1}$ and $A_{(n+1)}^{-1}$ are {\it row}  matrices of vectors. Taking the outer product of basis vectors in the relation 
$F_{(n+1)}=T^{-1}A_{(n+1)}$, gives
\[ \w F_{(n+1)}=\det T^{-1}\w A_{(n+1)}, \] 
or equivalently, after calculating and simplifying,
\beq e_1f_1 \cdots f_n = -\frac{(\sqrt 2)^{n+1}}{\sqrt n}a_1 \w \cdots \w a_{n+1}, \label{pseudoscalars} \eeq
relating the pseudoscalar elements of the geometric algebra $\G_{1,n}$ expressed in the standard basis and in the null vector basis of ${\cal A}^+_{n+1}$.
 
 \subsection{Decomposition formulas for $\nabla$}
 
 Usually the concept of {\it duality} is defined in terms of the operation of multiplication in an algebraic structure. By defining the geometric algebras $\G_{1,n}$ and $\G_{n,1}$ in terms of the null vector basis algebras ${\cal A}_{1,n}^+$ and ${\cal A}_{n,1}^-$, whose rules of multiplication have been given in the Multiplication Tables 1 and 2, suggests defining the concept of duality in terms of addition. The {\it dual $n$-sum} ${\vup a_i}$ of $a_i\in {\cal A}_{n+1}^+$ is the $n$-sum
 \beq {\vup a_i}:=  a_1+ \cdots + {\vup i} + \cdots a_{n+1} ,  \label{dualnsum} \eeq
 formed leaving out the $i^{th}$ term of the basis
 null vectors $\{a_1, \ldots, a_{n+1} \}\subset {\cal A}_{n+1}^+$. 

Calculations with the gradient $\nabla$ in ${\cal A}_{n+1}^+$ can often be simplified using the following {\it decompostion formulas}. Defining the {\it dual sum} and {\it null gradients} 
\beq  \vnabla := \sum_{i=1}^{n+1} {\vup a_i}\partial_i \ \ {\rm and} \ \  \hat \nabla := \sum_{i=1}^{n+1}  a_i\partial_i, \label{vnabla} \eeq
respectively, the gradient $\nabla$ defined by (\ref{gradientx}), $\vnabla$ and $\hat \nabla$ defined in (\ref{vnabla}), satisfy the following decomposition formulas:
\begin{itemize}
 \item $\nabla = \frac{2}{n}\Big(A_{n+1}\partial_{(n+1)}-n\hat \nabla\Big) =\frac{2}{n}\Big({\vnabla} - (n-1)\hat \nabla \Big)$, 
 
 where $\partial_{(n+1)} := \sum_{i=1}^{n+1} \partial_i$.
 \item $A_{n+1}\cdot \nabla =(n+1)\partial_{(n+1)} - 2 A_{n+1}\cdot {\hat \nabla}$
  \item $\vnabla + \hat \nabla = A_{n+1}\partial_{(n+1)}$  \ $\iff$  \ $A_{n+1}\cdot {\vnabla} + A_{n+1}\cdot \hat \nabla = \frac{(n+1)n}{2}\partial_{(n+1)}$

 \item $\hat\nabla^2 = \sum_{i<j}^{n+1}\partial_i \partial_j 	$, \ ${\vnabla}^2 = \frac{(n+1)n}{2}\sum_{i=1}^{n+1}\partial_i^2 + (n^2-n+1)\sum_{i<j}^{n+1} \partial_i \partial_j$
 
 \item $\nabla^2 =\frac{4}{n^2}\Big({\vnabla} - (n-1)\hat \nabla \Big)^2 = {\vnabla}^2 -2(n-1){\vnabla}\cdot \hat \nabla + \hat \nabla^2$,
 
   where ${\vnabla}\cdot \hat \nabla=\frac{n}{2}{\vnabla}\partial_{(n+1)}-\hat \nabla^2 $. 
\end{itemize}
Verifications of the above formulas, which are omitted, depend heavily on the combinatorial-like identities
\beq A_{n+1}^2 = \frac{(n+1)n}{2}, \ \ a_i \cdot A_{n+1} = A_{n+1} \cdot a_i = \frac{n}{2}{\vup a}_i , \label{Akident} \eeq  
and the {\it additive duality} formula for $A_{n+1}$, and $1\le i < j \le n+1$,
\beq {\vup a}_i \cdot {\vup a}_j = n^2 - n +1.  \label{Akident2} \eeq  

\section{Lower dimensional geometric algebras}
This section characterizes geometric sub-algebras of ${\cal A}_{3}^+\equiv \G_{1,2}$ in $\R^3$. 
 
 The pseudoscalar 
\beq i:=e_1f_1f_2 = -2a_1\w a_2 \w a_3 , \label{pseudoig3} \eeq
is in the center of the algebra, commuting with all elements. The algebra 
\[  \G_{3}:=\R(e_1,e_2,e_3), \]
is obtained from the algebra $\G_{1,2}$, simply by defining $e_2=i f_1=e_1f_2 \in \G_{1,2}^2$ and $e_3 = -i f_2 =e_1f_1\in \G_{1,2}^2$, and reinterpreting these anticommuting elements to be vectors in $\G_3^1$.

 The {\it matrix coordinates} $[e_1],[e_2],[e_3]$ of $e_1,e_2,e_3$, known as the famous {\it Pauli matrices}, opened the door to the study of quantum mechanics \cite[p.108]{GS2}. It has found many recent applications in computer science and robotics, \cite{Hitzer2023}. The geometric algebra $\G_{3}$ is isomorphic to the even subalgebra of the
{\it spacetime algebra} $\G_{1,3}={\cal A}_{4}^+ $ of $\R^4$. Its matrix version is known as the {\it Dirac algebra}. 

The null vector basis algebra ${\cal A}_3^+ =\G_{1,2}$
is defined by $3$ null vectors $\{a_1,a_2,a_3\}$ with the property that $\w A_{(3)}\ne 0$,  $a_i \cdot a_j = \frac{(1-\delta{ij})}{2}$, and the Multiplication Table 1. The relations between the standard basis of $\G_{1,2}$, and the basis of ${\cal A}_3^+$, are summarized by the $3\times 3$ transition matrix $T_3$, and its inverse,
\beq   T_3:=\pmatrix{\frac{1}{2} & \frac{1}{2} & 0 \cr \frac{1}{2} &- \frac{1}{2} & 0  \cr 1 & 0 & 1  },
\quad T_3^{-1}:=\pmatrix{1 & 1 & 0 \cr 1 &- 1 & 0  \cr -1 & -1 & 1  }. \label{transition4} \eeq

Using the relations given after (\ref{basesrel}), 
\beq \pmatrix{e_1 \cr f_1 \cr f_2 } = T_3^{-1} \pmatrix{a_1 \cr a_2 \cr a_3 } =\pmatrix{a_1+a_2 \cr a_1 - a_2 \cr -a_1 -a_2 +a_3 }  \label{stox3} \eeq
and 
\beq \pmatrix{a_1 \cr a_2 \cr a_3} = T_3 \pmatrix{e_1 \cr f_1\cr f_2 } =\pmatrix{\frac{1}{2}(e_1+f_1) \cr \frac{1}{2}(e_1-f_1)\cr e_1+f_2 }.    \label{xtos3} \eeq
The canonical forms relating the vectors, bivectors and trivectors are:
\begin{itemize}
	\item $e_1=a_1+a_2, \ f_1=a_1-f_1, f_2 = -a_1 -a_2+a_3$
	\item $e_1f_1=(a_1+a_2)(a_1-a_2)=a_2 a_1-a_1 a_2 = 1-2a_1a_2$
	\item $e_1f_2=(a_1+a_2)(-a_1-a_2+a_3)= 1+a_1a_3-a_2 a_3$
	\item $e_1f_1f_2=(a_1+a_2)(1-a_1a_2 + a_1a_3-a_2 a_3)=a_1+a_3-2a_1a_2a_3$ 
\end{itemize}

One of the simplest endomorphisms,
$f:\R^2 \to \R^{2}$,
defined by
$v_1,v_2 \in \R^2 $, is
\beq f(x):= 2(v_1\w v_2)x=2\Big((x\cdot v_2)v_1 - (x\cdot v_1)v_2\Big) , \label{deftoS2} \eeq
where $v_i = v_{i1}a_1 +v_{i2}$ for $i\in \{1,2\}$.
The endomorphism $f(x)$ has the {\it eigenvectors} $a_1$ and $a_2$, with the eigenvalues
$\pm \det\pmatrix{v_{11} & v_{12} \cr v_{12} & v_{22}} $,
\beq f(a_1)=2(v_1\w v_2)a_1 = \det \pmatrix{v_{11} & v_{12} \cr v_{12} & v_{22}} a_1 ,  \label{eigenvec1} \eeq  

\beq f(a_2)=2(v_1\w v_2)a_2 =-\det \pmatrix{v_{11} & v_{12} \cr v_{12} & v_{22}} a_2 ,  \label{eigenvec2} \eeq  
respectively, as is easily verified.

Now calculate,
\[ (v_1\w v_2)(v_1\w v_2 \w x ) =(v_1\w v_2)\cdot (v_1\w v_2 \w x )    \]
\[=(v_1\w v_2)^2 x+(v_1\w v_2)\cdot (v_2 \w x)v_1 +(v_1\w v_2)\cdot (x\w v_1)v_2 =0.         \]
Dividing both sides of this last equation by $(v_1\w v_2)^2$, gives
\beq \frac{(v_1\w v_2 \w x )}{(v_1\w v_2)}=x-\frac{(x \w v_2)}{(v_1\w v_2)}v_1 +\frac{(x\w v_1)}{(v_1\w v_2)}v_2 =0,   \label{xident} \eeq
expressing the position vector  $x \in \R^2$ uniquely in terms of its LPGG projective coordinates.
Of course, the trivector $v_1\w v_2 \w x=0$, because we are in the geometric algebra $\G_{1,1}$ of $\R^2$. 
Multiplying equation (\ref{xident}) by $4(v_1 \w v_2)^2$ immediately gives what I call the {\it Cayley-Grassmann identity},
\[ 4(v_1\w v_2)(v_1\w v_2 \w x )=f^2(x)-4(v_1\w v_2)(x \w v_2)v_1 +4(v_1\w v_2)(x\w v_1)v_2 \]
\beq =f^2(x)-2 f(x)\cdot v_2 v_1+ 2 f(x)\cdot v_1 v_2=0.   \label{CHident} \eeq

 The {\it matrix} of $[f(x)]$ of $f(x)$, in this translation, is given by
\beq [f(x)]=[2(v_1\w v_2)x]=2[(v_1\w v_2)][x],   \label{matrixfx} \eeq 
is the product of the matrix 
\beq [v_1 \w v_2]=\frac{1}{2}\Big([v_1v_2-v_2 v_1]\Big)=\frac{1}{2}\Big([v_1][v_2]-[v_2][ v_1]\Big)  ,\label{matrixbi} \eeq
where $[v_1\w v_2], [v_1], [v_2]$ are the matrices of $v_1\w v_2, v_1,v_2$, respectively,
and $[x]$ is the matrix of $x$. These matrices are given below.
With (\ref{stox3}) and (\ref{xtos3}) in hand, the matrix $[x]$ of the position vector $x\in \R^3$ 
\[ [x] =x_1[a_1]+x_2[a_2]+x_3 [a_3] = 
\pmatrix{x_3 i & x_2-x_3 \cr x_1 -x_3 & -x_3 i} \]
with respect to the basis ${\cal A}_3^+$, and
\[ [x] =s_1[e_1]+s_2[f_1]+s_3 [f_2] = 
\pmatrix{s_3 i & s_1 - s_2 \cr s_1 + s_2 & -s_3 i} \]
with respect to the standard basis of $\G_{1,2}$.

The $2\times 2$ matrices are defined with respect to the {\it spectral basis} 
\[ \pmatrix{a_2 a_1 & a_2 \cr a_1 & a_1 a_2}\]
of $\G_{1,1}$, as detailed in \cite{GS2} and \cite[p.78]{SNF}. The matrices of $[a_1]$ and $[a_2]$ of $a_1$ and $a_2$, are
\[ [a_1]=\pmatrix{0 & 0 \cr 1 &0} , \ \ [a_2]=\pmatrix{0 & 1 \cr 0 &0} ,   \]
respectively, and
\[ [x]=\pmatrix{0 & x_2 \cr x_1 &0}, \ [v_1]=\pmatrix{0 & v_{12} \cr v_{11} &0}, \ [v_2]=\pmatrix{0 & v_{22} \cr v_{21} &0},   \]
which are used with (\ref{matrixbi}) to calculate
\[ [v_1\w v_2]=\frac{1}{2}\Big([v_1][v_2]-[v_2][ v_1]\Big)= \frac{1}{2}\pmatrix{v_{12}v_{21}-v_{11}v_{22}&0 \cr 0 & v_{11}v_{22}-v_{12}v_{22}},   \]
and 
\[ [f(x)]=2[(v_1\w v_2)][x] = \pmatrix{0 &(v_{12}v_{21}-v_{11}v_{22})x_2 \cr -(v_{12}v_{21}-v_{11}v_{22})x_1   & 0}. \] 

Unlike the usual representation of an endomorphism $f(x)$ on $\R^2$, as a $2 \times 2$ matrix of $[f]$ of $f(x)$  times the column matrix of $x$, $[x]=\pmatrix{x_1 \cr x_2}$, the matrix of the endomorphism $f(x)$ in the LPGG of ${\cal V}_2^+(v_1,v_2)$ comes as the single real matrix $[f(x)]$. This single matrix $[f(x)]$ in LPGG can be broken into the product of two $2 \times 2$ matrices. By  (\ref{matrixfx}),
\[ [f(x)]= 2[(v_1\w v_2)][x]=\pmatrix{v_{12}v_{21}-v_{11}v_{22}&0 \cr 0 &-(v_{12}v_{21}-v_{11}v_{22}) }\pmatrix{0 & x_2 \cr x_1 & 0}. \]

For $k\in \{1,2,3\}$, let 
\[ v_k:=v_{k1}a_1+v_{k2}a_2+v_{k3}a_3 \in \R^3, \]
 consider the endomorphism
\beq f:\R^3  \to {\cal A}_3^+,\label{deftoS3} \eeq
defined by
\[ f(x):= 2(v_1 \w v_2 \w v_3 ) x = 2\det [v_{ij}] (a_1\w a_2 \w a_3 )x \]
\[=\det [v_{ij}]\Big((x_1+x_2 )a_1\w a_2+(x_2+x_3)a_2\w a_3+ (x_1+x_3) a_3\w a_1\Big).  \] 
It is interesting to note that each of the bivectors in the above expression are anticommutative and square to $\frac{1}{4}$. In view of (\ref{pseudoig3}), this is not surprising. Indeed, the mapping (\ref{deftoS3}) can simply be expressed as the duality relation $f(x)=-ix$. It follows that over the complex numbers, every vector $x\in {\cal A}_3^+$ is an eigenvector.

Consider the mapping $g:\R^3 \to {\cal A}_3^+$, defined by
\beq g(x):= \pmatrix{1&1&1}\pmatrix{0 &g_{12} a_1 a_2 & g_{13}a_1 a_3  \cr g_{21} a_2 a_1 & 0 &	g_{23}a_2 a_3 \cr  g_{31} 	a_{3}a_1 & g_{32}a_{3} a_2 & 0 }\pmatrix{1 \cr 1 \cr 1 }x=Gx,    \label{genbimap} \eeq 
where
\beq G=\frac{1}{2} tr(G)+{\rm g_1}a_2\w a_3+{\rm g_2}a_3\w a_1+{\rm g}_3a_1\w a_2,  \label{rep-pauli} \eeq
for ${\rm g}_1:=(g_{23}-g_{32}) $, ${\rm g}_2:= (g_{31}-g_{13})$, ${\rm g}_3:=(g_{12}-g_{21}) $, and
\[ tr(G):= g_{12}+g_{13}+g_{21}+ g_{23}+g_{31}+g_{32}.  \]
The same mapping (\ref{genbimap}) can equally well be considered over $\C$, 
\beq g:\C ^3 \to {\cal A}_3^+(\C), \label{paulig3} \eeq 
giving a new relationship between the Pauli matrices and $\G_3$. In this case, 
\[ e_1:=\frac{1}{2}a_2 \w a_3, \ e_2:=\frac{1}{2}a_3 \w a_1, \ e_3:=\frac{1}{2}a_1 \w a_2.   \]

The {\it minimal polynomial} of $G$ is easily calculated:
\[\varphi(G)= \big(G-\frac{1}{2}tr(G)\big)^2 - \frac{1}{4}\Big(  {\rm g_1}^2+{\rm g_2}^2+{\rm g_3}^2 \Big). \]
Setting $\varphi(r)=0$ and solving for $r$, gives the two roots $r_-$ and $r_+$,
 \[r_\mp:=\frac{1}{2}\Big(tr(G)\mp \sqrt{{\rm g}_1^2+{\rm g}_2^2+{\rm g}_3^2} \Big). \]          
In the {\it spectral basis}, $G$ takes the form
\beq G = r_- p_1(G) + r_+ p_2(G), \label{p1andp2} \eeq    
where
\[p_1(t):= \frac{-2t+tr(G)+\sqrt{{\rm g}_1^2+{\rm g}_2^2+{\rm g}_3^2}}{2\sqrt{{\rm g}_1^2+{\rm g}_2^2+{\rm g}_3^2}} \] 
and 
\[ p_2(t):= \frac{2t-tr(G)+\sqrt{{\rm g}_1^2+{\rm g}_2^2+{\rm g}_3^2}}{2\sqrt{{\rm g}_1^2+{\rm g}_2^2+{\rm g}_3^2}},  \]
\cite{SNF,GS2,college}.

\section{Simplices  in ${\cal A}^+_{n+1}$   }

It has been shown in previous sections how the  development of linear algebra can be carried out in $\R^{n+1}$, using the tools of $\G_{1,n}\equiv {\cal A}_{n+1}^+$. Restricting to barycentric coordinates, gives new tools for application in graph theory.  
In {\it Simplicial Calculus with Geometric Algebra}, many ideas of simplicial geometry were set down in the context of geometric algebra \cite{simcal92}. The present work is in many ways a continuation of this earlier work.

Let ${\cal A}_{n+1}^+$ be the null vector algebra of the geometric algebra $\G_{1,n}$, defined by the Multiplication Table 1, and where the null vectors $a_{i}$ satisfy for $1 \le i, j \le n+1$,
\beq a_i \cdot a_j = \frac{(1-\delta_{ij})}{2}. \label{defseqs} \eeq
For $x\in \R^{n+1}$ the position vector (\ref{xinA1n}), the {\it convex null}
{\it $n$-simplex} in $\R^{n+1}$ is defined by
\beq {\cal S}^+_n:={\cal S}_n^+(a_1, \ldots, a_{n+1})=\{x\in \R^{n+1}| \ x_1+\cdots +x_{n+1} = 1,  \ x_i\ge 0\}, \label{defSn} \eeq
by the requirement that the coordinates $x_{(s)}$ of $x\in \R^{n+1}$, are {\it homogeneous barycentric coordinates}, \cite{wikibary}.  

By the {\it content} of ${\cal S}^+_n$, we mean
\[  a_{{\triangle}_n}:=\frac{1}{n!} \w_{i=2}^{(n+1)}(a_i-a_1)= (a_2-a_1)\w (a_3-a_1)\w   \cdots \w (a_{n+1}-a_1) \]
\beq =\frac{1}{n!}\Big( \w {\vup a_1} - \w {\vup a_2} + \cdots +(-1)^{n} \w {\vup a_{n+1}} \Big).  \label{simcontent} \eeq 
Wedging (\ref{simcontent}) on the left by $x\in  {\cal S}_{n}^+$, gives
\beq x \w a_{{\triangle}_n} =\frac{1}{n!}\Big(\sum_{i=1}^{n+1}x_i\Big) \w A_{(n+1)} =\frac{1}{n!} \w A_{(n+1)}.  \label{altbary} \eeq
Similarly, dotting (\ref{simcontent}) on the left by $x$ gives
\beq x \cdot a_{{\triangle}_n} =\frac{1}{n!}x \cdot \Big( \w {\vup a_1} -\w {\vup a_2} + \cdots +(-1)^{n} \w {\vup a_{n+1}} \Big).   \label{dotbary} \eeq

Let $v_1, \ldots, v_{k+1} \in \R^{k+1}$ be a set of $k+1$ vertices of a {\it $k$-simplex} ${\cal V}_k^{+}$ in ${\cal A}_{k+1}^+$. That is 
\beq v_i := \sum_{j=1}^{k+1} v_{ij} a_j = [v_{ij}]A_{(k+1)} ,\label{simplexpts} \eeq
where $[V]_{k+1}:=[v_{ij}]$ is the {\it matrix} of ${\cal V}_k^{+}$. The {\it rows} of the {\it simplicial matrix} $[V]_{k+1}$ are the barycentric coordinates of the vertices $v_i \in {\cal V}^+_{k+1}$. It follows that $[V]_{k+1}$ is a non-negative matrix with the property that the sum of the coordinates in each row is equal to $1$. Alternatively, since $v_1 \w \cdots \w v_{k+1} \ne 0$, and not requiring the coordinates to be barycentric, the matrix $[V]_{k+1}$ becomes the transition matrix from the basis of null vectors ${\cal A}_{k+1}$ to 
the basis vectors $v_i \in {\cal V}_{k+1}^+$, for which all the relations found after (\ref{basesrel}) remain valid. 

 The {\it content} of ${\cal V}_{k+1}^+$ is
\beq   v_{{\triangle}_k} =   \w_{i=2}^{k+1}(v_i-v_1)= (v_2-v_1)\w  \cdots \w (v_{k+1}-v_1) \ne 0 ,  \label{regkpoly} \eeq
in the geometric algebra ${\cal A}_{k+1}^+$ of $\R^{k+1}$. Similar to (\ref{altbary}) and (\ref{dotbary}), we have
\[ xv_{{\triangle}_k}=x \cdot v_{{\triangle}_k} + x\w  v_{{\triangle}_k},  \]
but there is no obvious simplification as found for null simplices in (\ref{altbary}).

\subsection{LPGG Calculus of ${S}^+_{n}\subset \R^{n+1}$} 

I will now give a brief introduction to the general theory of {\it LPGG Calculus}. Standard {\it geometric calculus} has been in continual development over the last half Century \cite{HS1992,SNF,PL2001,Hitzer2023}.
Every {\it signed graph} ${\cal V}_n^\pm $ of $n$-vertices can be studied in terms of any of the geometric algebras determined by the sequences of signs
$(\ref{sequence1})$, (\ref{sequence2}), found in Appendix B. I will limit my discussion here to {\it signed positive $\frac{1}{2}$-graphs} ${\cal V}_m^+$ in $\R^{n+1}$, using the barycentric coordinates of the convex null symplex  
\[S_{n}^+:={\cal S}_{n}^+(a_1,\ldots, a_{n+1})\subset {\cal N},\]
  and the geometric algebra ${\cal A}_{n+1}^+\equiv \G_{1,n}$. For $n=0$, define 
\[ S_0^+:=\{ a |\ a \ne 0, a^2=0 \} \subset {\cal N} \]
to be the single null vector $a$.

  For $m\le n$, let $v_1, \ldots v_{m+1} \in \R^{n+1}$ 
denote the vertices of a signed simplex 
\[ {\cal V}_{m}^+:={\cal V}_{m}^+(v_1,\ldots, v_{m+1})\subset {\cal S}_n^+,\]
 where $v_1 \w \cdots \w v_{m+1} \ne 0$. 
 Define 
 $a_{(n+1)}\equiv\{a \}_{(n+1)}$ by
\[\{a \}_{(n+1)}:=\{a_1,\ldots,a_{n+1}\},   \]
and by $\{\vup a_i \}_{(n)}$, the set of $n$ correlated null vectors obtained by leaving out $a_i$, 
\[ \{\vup a_i \}_{(n)} :=\{a_1,\ldots \vup a_{i} \ldots, a_{n+1}\}. \]
When no confusion can arise, we shorten $\{a\}_{(n)}$ to
$a_{(n)}$. 
For $n=3$,
\[ \{\vup a_1 \}_{(3)}=\{a_2, a_3\},\ \{\vup a_2 \}_{(3)} =\{a_1, a_3\},\ {\rm and} \  \{\vup a_3 \}_{(3)}=\{a_1, a_2\}.  \]

Since the set of vectors $ \{v \}_{(m+1)}$ are {\it linearly independent}, 
\[ \w v_{(m+1)}:= v_1 \w \cdots \w v_{m+1} \ne 0 ,\]
 ${\cal V}_{m}^+$ defines an 
$m$-simplex with $m+1= \pmatrix{m+1 \cr m}  
$-{\it faces}. Each $(m+1)$-face is geometrically  represented by the oriented $m$-vector  
\[ \w\vup v_{(i)}:=v_1\w \cdots \vup i  \cdots   \w v_{m+1}. \]
Also, define the $(m+1)$-sum and the $m$-sum, by
\[ \sum v_{(m+1)}:=v_1 + \cdots + v_{m+1}, \ {\rm and} \ \sum v_{(\vup i)}:=v_1 + \cdots \vup i  \cdots + v_{m+1}.  \]

The {\it signed complete graph} ${\cal V}_{m}^+$ is said to be {\it closed} if
$ \sum_i v_{(\vup i)}=0$, and of {\it order} $k$, if
$k$ is the largest number of linearly independent
vertices of ${\cal V}^+_{m}$. 
Naturally, we use the {\it barycentric coordinates} associated with ${\cal S}^+_{n}$, and, without loss of generality, assume that the {\it position vector} $x \in {\cal V}^+_{m+1}$, is given by
\[ x=(x_1, \cdots ,x_{m+1}):= \sum_{i=1}^{m+1} x_i a_i \in \R^{m+1} \subset \R^{n+1} ,\]
 although other coordinate systems can be used.

I now restrict attention to studying the graph ${\cal V}^+_{m+1}$ of a particular $m$-dimensional {\it polytope}. For $x \in {\cal V}^+_{m+1}$, calculate
\[ x^2=\bigg( \sum_{i=1}^{m+1} x_i a_i \bigg)^2=\sum_{0\le i<j \le m+1} x_i x_j,  \]
for all $i,j$, $0<i\ne j \le m+1$, and where the vertices of the $m$-polytope satisfy
\[  v_1 \w \cdots \w v_{m+1} \ne 0. \]

Since ${\cal V}_m^+ \subset {\cal S}_n^+$, the barycentric coordinates $x_i$ of $x$ will all be positive, so that
\beq |x|^2=x^2 = \sum_{0\le i<j \le m+1} x_i x_j \ge 0. \label{squarexbary} \eeq   
For $x \in {{\cal V}}^+_{m}$, define
\[ |x| =\sqrt{ \sum_{0\le i<j \le m+1} x_i x_j}\ \ge 0. \] 
The points $x \in {\cal V}^+_{m}\subset {\cal S}_n^+$, for which $|x|=0$, are exactly those points $x$ of the graph on the light cone. 
For all interior points of ${\cal V}_m^+$, where $|x|>0$,
define the unit vector
\beq \hat x := \frac{x}{|x|}.   \label{Aunitvec} \eeq

Since $x\in {\cal S}_n^+$ is barycentric, its coordinates satisfy $\sum x_i = 1$. Taking the partial derivative $\partial_i$ of this equation, gives $\partial_i \sum_{j=1}^{n+1} x_j =0$. By employing {\it higher order barycentric coordinates}, based upon {\it Hermite interpolation}, this constraint can be satisfied. Without going into details, for $x\in {\cal S}_n^+$, I want to preserve the property that $\partial_i x=a_i$ for each $0<i \le n+1$, \cite{college,highbary}. The same effect can be achieved by assuming, when differentiating $x$, we have relaxed the condition that the coordinates of $x\in \R^{n+1}$ are barycentric. 

Recalling (\ref{gradientx}) and (\ref{vnabla})
\beq \nabla = \frac{2}{n}\Big(A_{n+1}\partial_{(n+1)}-n\hat \nabla\Big) =\frac{2}{n}\Big({\vnabla} - (n-1)\hat \nabla \Big). \label{nabladecomp} \eeq
For $x\in \R^{n+1}, \nabla x^2 = 2x$, $\nabla |x|=\hat x$, $\nabla |x|=\hat x$, and $\nabla \hat x=\frac{n}{|x|}$. These formulas remain valid at all points $x\in {\cal S}_{n}^+$, \cite[p.66]{SNF}. 
\subsection{Platonic solids}

Applying the decomposition formula (\ref{nabladecomp}),
\[ \vnabla x =  \sum_i {\vup a}_{i}  \frac{\partial x}{\partial i} =  \sum_i {\vup a}_{i} a_i=  \pmatrix{n\cr 2}= \frac{n(n-1)}{2},\]
which is the number of linear independent edges of ${S}^+_{n}$. 

The {\it Laplacian} $\vnabla^2$ for the light cone projective geometry of ${\cal S}^+_{n}$ is 
\[ \vnabla^2 =\sum_{i=1}^n  \partial_i^2+  \pmatrix{n \cr 2}\sum_{1\le i \le j\le n}\partial_i \partial_j . \] 
Just as in Euclidean and pseudo-Euclidean geometry,
the Laplacian $\vnabla^2$ in ${\cal S}^+_{n}$, is scalar valued.

For the signed simplex ${\cal S}^+_{3}$,
\[ \vnabla^2 =  \big(\partial_1^2+\partial_2^2+\partial_3^2 +\partial_2\partial_3+\partial_1\partial_3+\partial_1 \partial_2\big).\]
For $x\in {\cal S}^+_{n} $, we calculate
\[ \vnabla^2 x^2 = \sum_{i=1}^n  \partial_i^2x^2+ \pmatrix{n \cr 2}\sum_{1\le i \le j\le n}\partial_i \partial_j x^2=\pmatrix{n \cr 2}^2. \]

I conclude with a {\it Conjecture} for {\it $n$-Platonic Solids} in $(n+1)$-dimensional space $\R^{n+1}$.

 {\bf Conjecture:}
  {\it The number of $n$-Platonic Solids in any  dimension $n$ is equal to the number of distinct $n$-Platonic Solids found in the $n$-simplex ${\cal S}^+_{n}\subset \R^{n+1}$ with its vertices located at the null vectors $a_1, \ldots, a_{n+1} \in {\cal S}^+_{n}$.}

  The number is known to be given by the sequence
	\[\{ 1,1,\infty, 5,6,3,3,3, \cdots \},  \]
	\cite{LPGG-Nov16,Havel2022,Tanya2008,baez,Hest2006}.

I want to welcome the reader to this beautiful new, but not really so new, theory. Be careful - the calculations can be treacherous.

\section*{Acknowledgements} 
The seeds of this note were planted almost 40 years ago in discussions with
 Professor Zbigniew Oziewicz, a distinguished colleague,
 about the fundamental role played by {\it duality} in its many different guises in mathematics and physics \cite{ZO1985}.  The author thanks the {\it Zbigniew Oziewicz Seminar on Fundamental Problems in Physics} group for many fruitful discussions of the ideas herein \cite{FES-C}, and offers special thanks to Timothy Havel for thoughtful comments about earlier versions of this work.
Not least, the author thanks the organizers of ICACGA2023 for allowing me extra time to complete this work.

\section*{Appendix A: Geometric Algebra Identities in ${\cal A}_{1,n}^+$}

Some basis identities of the geometric algebra 
\[  \G_{1,n}\equiv {\cal A}_{1,n}^+=\R^{n+1}:=\R(a_1, \ldots , a_{n+1}),\]
 where $a_{i}\cdot a_j = \frac{1-\delta_{ij}}{2}$.

\begin{itemize}
	\item[1.] $x^2=x_1 x_2, \  x\cdot v_1=\frac{1}{2} (x_1 v_{12}+x_2 v_{11} ), \  x\cdot v_2=\frac{1}{2} (x_2 v_{22}+x_2 v_{21} )$
	\item[2.] $v_1\cdot v_2=\frac{1}{2} (v_{11} v_{22}+v_{12} v_{21} ), \ v_1\w v_2=\det \pmatrix{v_{11}& v_{12} \cr v_{21} & v_{22}}a_1\w a_2$
	\item[3.] $(a_1\w a_2)=\frac{1}{2}(a_1-a_1)\w (a_1+a_2)=\frac{1}{2}f_1e_1, \ (a_1\w a_2)^2=\frac{1}{4}, $
	\item[4.] For $y=y_1 a_1+y_2 a_2$, $ x\w y= \det\pmatrix{x_1 & x_2 \cr y_1 & y_2}a_1\w a_2$
	\item[5.] $(x\w y)^2=\det\pmatrix{y\cdot x & y^2 \cr x^2 & x\cdot y}$
\end{itemize}  

\centerline{\large\bf Change of Basis Formulas for $\bf n+1=8$}

\[ T_8=  \pmatrix{\frac{1}{2} &\frac{1}{2}  & 0 &0 &0&0 &0 &0  \cr 
	\frac{1}{2} &-\frac{1}{2}  & 0 &0 &0&0 &0 &0 \cr 1 & 0  & 1 &0 &0&0 &0 &0 \cr 1 & 0  & \frac{1}{2} &\frac{\sqrt 3}{2} &0&0 &0 &0 \cr
	 1 & 0  & \frac{1}{2} &\frac{1}{2 \sqrt 3}&\sqrt{\frac{ 2}{3}} &0 &0&0  
	 \cr  1 & 0  & \frac{1}{2} &\frac{1}{2 \sqrt 3}&\frac{ 1}{2\sqrt{6}} &\frac{\sqrt 5}{2\sqrt 2} &0&0 
	 \cr  1 & 0  & \frac{1}{2} &\frac{1}{2 \sqrt 3}&\frac{ 1}{2\sqrt{6}} &\frac{1}{2\sqrt 10} &\sqrt{\frac{ 3}{ 5}}&0 
	 \cr   1 & 0  & \frac{1}{2} &\frac{1}{2 \sqrt 3}&\frac{ 1}{2\sqrt{6}} &\frac{1}{2\sqrt 10} &\frac{ 1}{2\sqrt{15}}&\frac{\sqrt 7}{2\sqrt 3} } \]
 
 \[ T_8^{-1}=  \pmatrix{1 &1  & 0 &0 &0&0 &0 &0  \cr 
 	1 &-1  & 0 &0 &0&0 &0 &0 
 	\cr -1 & -1  & 1 &0 &0&0 &0 &0 
 	\cr -{\frac{ 1}{\sqrt 3}} &-{\frac{ 1}{\sqrt 3}}  &-{\frac{ 1}{\sqrt 3}}&-{\frac{ 2}{\sqrt 3}} &0&0 &0 &0
 	 \cr
 -{\frac{ 1}{\sqrt 6}} &-{\frac{ 1}{\sqrt 6}}  &-{\frac{ 1}{\sqrt 6}}&-{\frac{ 1}{\sqrt 6}} &\sqrt \frac{3}{2}&0 &0 &0
 	\cr   -{\frac{ 1}{\sqrt 10}} &-{\frac{ 1}{\sqrt 10}}  &-{\frac{ 1}{\sqrt 10}}&-{\frac{ 1}{\sqrt 10}} &-\frac{1}{\sqrt 10}&2 \sqrt \frac{2}{5} &0 &0
 	\cr -{\frac{ 1}{\sqrt 15}} &-{\frac{ 1}{\sqrt 15}}  &-{\frac{ 1}{\sqrt 15}}&-{\frac{ 1}{\sqrt 15}} &-\frac{1}{\sqrt 15}&-\frac{1}{\sqrt 15} &\sqrt \frac{5}{3} &0
 	\cr    -{\frac{ 1}{\sqrt 21}} &-{\frac{ 1}{\sqrt 21}}  &-{\frac{ 1}{\sqrt 21}}&-{\frac{ 1}{\sqrt 21}} &- \frac{1}{\sqrt 21}&- \frac{1}{\sqrt 21} &- \frac{1}{\sqrt 21} &2 \sqrt \frac{3}{7} } \]

\section*{Appendix B: Classification of Geometric Algebras}
There is an extremely interesting relationship between plus and minus signs of the squares of the standard basis elements of $\G_{p,q}$, and the $8$-fold periodicity structure of Clifford geometric algebras.
Consider the following:
\begin{itemize}
	
	\item[1.] $\{+ \}, \{-\},\ \ e_1 \in \G_{1,0}, \ f_1 \in \G_{0,1}  \hfill  \prod {signs}\  -$
	\item[2.] $\{+ +\},\  \{+ -\},\ \{--\},\ \ \G_{p,q}, \ \ p+q=2  \hfill  \prod {signs} \  -  $ 
	\item[3.] $\{+++\},\  \{++-\},\ \{+--\},\  \{---\} \ \ p+q=3, \ etc. \hfill   \prod {signs} \   +  $
	\item[4.] $\{++++\},\  \{+++-\},\ \{++--\},\  \{+---\},\  \{----\}\hfill  \prod {signs} \ +  $
	\item[5.] $\pmatrix{n+1 \cr 2}= \pmatrix{6 \cr 2}=15   \hfill \prod \ {(15)} \ -$	
	\item[6.] $ \pmatrix{n+1 \cr 2}= \pmatrix{7 \cr 2}=21   \hfill \prod {\ (21)} \ -$	
\end{itemize}
This obviously gives the infinite sequence
\beq  --,++,--,++,--,++\ldots .\label{sequence1} \eeq
{\it Real geometric algebras} $\G_{p,q}$ are constructed by extending the real number system $\R$ by $n=p+q$ anti-commuting vectors $e_i,f_j$ which have
sqares $\pm 1$, respectively
\beq \G_{p,q}:=\R[e_1,\cdots,e_p,f_1,\cdots f_q], \label{stdbasisGpq} \eeq
\cite{HS1992,SNF}. A more concise treatment of this construction, and its relationship to real and complex square matrices is \cite{GS2}.

Geometric algebras enjoy a very special $8$-fold periodicity relationship \cite{GS3}. A basic understanding of this important periodicity relationship can be obtained by studying the signs of the squares of the pseudoscalar elements for the geometric algebra $\G_{p,q}$ of successively higher dimensions. The $\pm$ signs over {\it pseudoscalar elements} indicate the sign of the {\it square} of that element.

\begin{itemize}
	\item[0.] $\{ a\in {\cal N}| \ a^2=0 \}$. The {\it null vector $a\ne 0$} has the property that $a^2=0$. 
	\item[1.] $\{ \stackrel{+}{e}_1 \}, \ \{\stackrel{-}f_1\}: \ \ \G_{1,0}, \ \G_{0,1};    $
	\item[2.] $\{\stackrel{-}{e_1e_2}\}\ \{\stackrel{+}{e_1f_1} \},\ \{\stackrel{-}{f_1f_2} \}: \ \G_{2,0}, \ \G_{1,1}, \ \G_{0,2}; $	
	\item[3.] $\{\stackrel{-}{e_1e_2e_3}\},\ \{\stackrel{+}{e_1e_2f_1} \}, \ \{\stackrel{-}{e_1f_1f_2} \}, \ \{\stackrel{+}{f_1f_2f_3} \}:\ \G_{3,0}, \ \G_{2,1}, \ \G_{1,2}, \ \G_{0,3};  $
	\item[4.] $\{\stackrel{+}{e_1e_2e_3e_4}\},\ \{\stackrel{-}{e_1e_2e_3f_1} \}, \ \{\stackrel{+}{e_1e_2f_1f_2} \}, \ \{\stackrel{-}{e_1f_1f_2f_3} \}, \ \{\stackrel{+}{f_1f_2f_3f_4} \} $
	\item[5.] $\{\stackrel{+}{e_1e_2e_3e_4e_5}\},\ \{\stackrel{-}{e_1e_2e_3e_4f_1} \}, \ \{\stackrel{+}{e_1e_2e_3f_1f_2} \}, \ \{\stackrel{-}{e_1e_2f_1f_2f_3} \}, $
	
	$ \{\stackrel{+}{e_1f_1f_2f_3f_4} \},\{\stackrel{-}{f_1f_2f_3f_4f_5} \}$.
	\item[6.] $\{\stackrel{-}{e_1e_2e_3e_4e_5e_6}\},\ \{\stackrel{+}{e_1e_2e_3e_4e_5f_1} \}, \ \{\stackrel{-}{e_1e_2e_3e_4f_1f_2} \}, \ \{\stackrel{+}{e_1e_2e_3f_1f_2f_3} \}, $
	
	$ \{\stackrel{-}{e_1e_2f_1f_2f_3f_4} \},\{\stackrel{+}{e_1f_1f_2f_3f_4f_5},\{\stackrel{-}{f_1f_1f_2f_3f_4f_5} \}$.
\end{itemize}
This obviously gives the sequence,
\beq +,-,-+-,-+-+,+-+-+,-+-+-+, \cdots  \label{sequence2} \eeq

The sequences (\ref{sequence1}) and (\ref{sequence2}) follow directly from the well known periodicity laws of all real and complex geometric algebras \cite{PL2001}. The two sequences beautifully reflect how any geometric algebra $\G_{p,q}$, for $n=p+q$ can be represented either as a real or complex matrix algebra of dimension $2^n$. In the case of the complex matrix algebra, the imaginary number $i$ can be interpreted as the pseudoscalar element
$e_1f_1 \cdots e_nf_n f_{n+1}$ in the center of the real geometric algebra $\G_{n,n+1}$.

\end{document}